\documentclass[11pt]{article}

\usepackage{amsthm,amsmath,amssymb,tikz,graphicx,array,arydshln,geometry}
\usetikzlibrary{patterns}
\usepackage[colorlinks=true,
            linkcolor=blue,
            urlcolor=blue,
            citecolor=blue]{hyperref}

\newtheorem{thm}{Theorem}
\newtheorem{lem}[thm]{Lemma}
\newtheorem{cor}[thm]{Corollary}

\theoremstyle{definition}

\newtheorem{pr}[thm]{Problem}

\theoremstyle{remark}
\newtheorem{rem}[thm]{Remark}
\newtheorem*{rem*}{Remark}

\tikzset{>=latex}

\newcommand{\floor}[1]{\left\lfloor #1\right\rfloor}
\newcommand{\ceil}[1]{\left\lceil #1\right\rceil}

\def\R{\mathbb{R}}

\def\le{\leqslant}
\def\ge{\geqslant}

\def\A#1{\mathcal{A}(#1)}

\def\bleq{\preceq_{B}}
\def\sble{\prec_{\widehat{B}}}
\def\sbleq{\preceq_{\widehat{B}}}
\def\cnj#1{\overline{#1}}
\def\w{w}
\def\h{h}

\def\etal{~{\it et al.}}

\begin{document}

\title{Antichains of $(0,1)$--matrices through inversions}

\author{M. Ghebleh\\[1ex]
{\em Department of Mathematics,
             Faculty of Science,
             Kuwait University,
             Kuwait}\\
\tt mohammad.ghebleh@ku.edu.kw}

\date{13 January 2014}

\maketitle

\begin{abstract}
An inversion in a matrix of zeros and ones consists of two entries
both of which equal $1$, and one of which is located to the top-right
of the other.
It is known that in the class $\mathcal{A}(R,S)$ of $(0,1)$--matrices with
row sum vector $R$ and column sum vector $S$, the number of inversions in a
matrix is monotonic with respect to the secondary Bruhat order.
Hence any two matrices in the same class $\mathcal{A}(R,S)$
having the same number of inversions, are incomparable in the
secondary Bruhat order.
We use this fact to construct antichains in the Bruhat order
of $\mathcal{A}(n,2)$, the class of all $n\times n$ binary matrices 
with common row and column sum~$2$.

A product construction of antichains in the Bruhat order of
$\mathcal{A}(R,S)$ is given. This product construction is applied
in finding antichains in the Bruhat order of the class $\mathcal{A}(2k,k)$
of square $(0,1)$--matrices of order $2k$ and common row and column sum~$k$.
\end{abstract}


\section{Introduction}

Let $R=\langle r_1,\ldots,r_m\rangle$ and $S=\langle s_1,\ldots,s_n\rangle$
be vectors of nonnegative integers.
The class of all $m$ by $n$ matrices of zeros and ones with
row sum vector $R$ and column sum vector $S$ is denoted by $\A{R,S}$.
The class of all square $(0,1)$--matrices of order $n$ with common
row and column sum $k$ is denoted by $\A{n,k}$.
Combinatorial properties of the class $\A{R,S}$ are studied
extensively (see for example \cite{Brualdi1980,Brualdi2006,BrualdiBook2006,GreenhillMcKayWang2006,Ryser1963}
and the references therein). An asymptotic formula for the
size of the class $\A{n,k}$ due to O'Neil~\cite{oneil} is reported
in~\cite{EverettStein}:
\begin{equation}
\left|\A{n,k}\right|\sim\frac{(kn)!}{(k!)^{2n}}\,e^{-(k-1)^2/2}
\label{eq:asymptotic}
\end{equation}

Brualdi and Hwang~\cite{BrualdiHwang2004} define a {\em Bruhat order}
on the class $\A{R,S}$ generalizing the classical Bruhat order
on the symmetric group $\mathcal{S}_n$.
To any $(0,1)$--matrix $A$ of size $m$ by $n$,
an $m$ by $n$ matrix $\Sigma_A$ is assigned whose $(i,j)$--entry is
\[
\sigma_{i j}(A)=\sum_{k=1}^{i}\sum_{\ell=1}^{j}a_{k\ell}.
\]
If $A$ and $C$ are $(0,1)$--matrices in a class $\A{R,S}$, then 
$A$ precedes $C$ in the Bruhat order, written $A\bleq C$,
if $\sigma_{ij}(A)\ge\sigma_{ij}(C)$ for all $1\le i\le m$ and 
$1\le j\le n$.
It is easy to see that if $C$ is obtained from $A$ by a sequence of
\[
I_2=\left[\begin{array}{cc}1&0\\0&1\end{array}\right]
\to
\left[\begin{array}{cc}0&1\\1&0\end{array}\right]=L_2
\]
submatrix interchanges, then $A\bleq C$. 
This defines a {\em secondary Bruhat order} on the class
$\A{R,S}$:
$A\sbleq C$ if and only if $C$ is obtained from $A$ by a sequence of
$I_2\to L_2$ interchanges.
It is shown in~\cite{BrualdiDeaett2007} that the Bruhat order and the
secondary Bruhat order are the same on the classes $\A{n,2}$,
but they are different on $\A{6,3}$.

Two classical theorems of Dilworth~\cite{dilworth} and Mirsky~\cite{mirsky}
give relationships between chains and antichains in a partially ordered set.
Dilworth's theorem states that the maximum cardinality of an antichain
in a partially ordered set equals the minimum number of chains into
which the set may be partitioned.
Mirsky's dual of this theorem states that the maximum cardinality
of a chain in a partially ordered set equals the minimum number of
antichains into which the set may be partitioned.
The \emph{height} (resp. \emph{width}) of a partially ordered set
$(P,\preceq)$ is defined to be the maximum cardinality of a chain
(resp. antichain) in~$(P,\preceq)$.
If $h$ and $w$ denote the height and the width of a finite poset
$(P,\preceq)$ respectively, Dilworth's and Mirsky's theorems imply
$hw\ge|P|$. Therefore, an upper bound on one of $h$ and $w$ gives
a lower bound on the other. In particular, if $\h(n,k)$ and $\w(n,k)$
denote the height and the width of the Bruhat order of $\A{n,k}$
respectively, then
\begin{equation}
\w(n,k)\ge\frac{\left|\A{n,k}\right|}{\h(n,k)}.
\label{eq:width}
\end{equation}
Conflitti\etal~\cite{ConflittiDaFonsecaMamede2012}
show that the Bruhat order of the class $\A{2k,k}$ where $k\ge1$
is an integer, has height $\h(2k,k)=k^4+1$.
In~\cite{Mamad:An2Chain} the height of the Bruhat order of $\A{n,2}$
for an integer $n\ge4$, is proved to be $\h(n,2)=2n(n-2)+\varepsilon$
where $\varepsilon$ equals $0$ if $n$ is odd, and $1$ if $n$ is even.
These polynomial values together with equations~(\ref{eq:asymptotic})
and~(\ref{eq:width}) indicate that the widths $\w(n,2)$ and $\w(2k,k)$
have exponentially large sizes.
Conflitti\etal~\cite{ConflittiDaFonsecaMamede2013} prove a
lower bound of order $k^8$ for the width of Bruhat order of $\A{2k,k}$.
In this work we show that
\[
\w(n,2)\ge f(n)=\begin{cases}
\displaystyle\frac{n!}{2^{n/2}} & \text{ if $n$ is even,}\\[1em]
\displaystyle\frac{(n-1)!}{2^{(n-3)/2}} & \text{ if $n$ is odd,}
\end{cases}
\]
and
\[
\w(2k,k)\ge g(k)=\begin{cases}
\displaystyle\frac{(k!)^4}{4^{k}} & \text{ if $k$ is even,}\\[1em]
\displaystyle\frac{\big[(k-1)!\big]^4}{4^{k-3}} & \text{ if $k$ is odd.}
\end{cases}
\]

\section{Antichains in the Bruhat order of $\A{n,2}$\label{sec:An2}}

The notion of an \emph{inversion} in a $(0,1)$--matrix $A=[a_{ij}]$
is introduced in~\cite{Mamad:An2Chain}.
Two entries $a_{ij}$ and $a_{k\ell}$ constitute an inversion in $A$,
if $a_{ij}=a_{k\ell}=1$, and $(i-k)(j-\ell)<0$.
The total number of inversions in $A$ is denoted by $\nu(A)$.
It is proved in~\cite{Mamad:An2Chain} that the number of inversions
is monotonic with respect to the secondary Bruhat order of any class
$\A{R,S}$. Namely, if $A,C\in\A{R,S}$ such that $A\sble C$, then
$\nu(A)<\nu(C)$. This implies that if $\nu(A)=\nu(C)$ for some
$A,C\in\A{R,S}$, then $A$ and $C$ are incomparable in the
secondary Bruhat order of $\A{R,S}$.
On the other hand, it is proved in~\cite{BrualdiHwang2004}
that the Bruhat order and the secondary Bruhat order are the same
on any class $\A{n,2}$. Thus we obtain the following.

\begin{lem}
Let $n\ge2$ and $t\ge0$ be integers.
Then the set $\nu^{-1}(t)$ of all matrices $A\in\A{n,2}$ with
$\nu(A)=t$ is an antichain in the Bruhat order of $\A{n,2}$.
\label{lem:inv_antichain}
\end{lem}

Note that depending on the value of $t$, the antichain of the
above lemma may be empty. It is proved in~\cite{Mamad:An2Chain} that
$\ceil{n/2}\le \nu(A)\le 2n(n-2)+\floor{n/2}$ for every $A\in\A{n,2}$,
and it is implicit there that for $n\ge4$ and $\ceil{n/2}\le t\le 2n(n-2)+\floor{n/2}$, the antichain $\nu^{-1}(t)$ of the above lemma
is nonempty.
The following is immediate from the above discussion.

\begin{cor}
Let $n\ge3$ be an integer.
Then $\displaystyle\w(n,2)\ge\frac{|\A{n,2}|}{2n(n-2)+\varepsilon}$,
where $\varepsilon$ equals $0$ if $n$ is odd and $1$ if $n$ is even.
\label{cor:average}
\end{cor}

\begin{figure}
\begin{center}
\begin{tabular}{c}
\includegraphics[scale=.58]{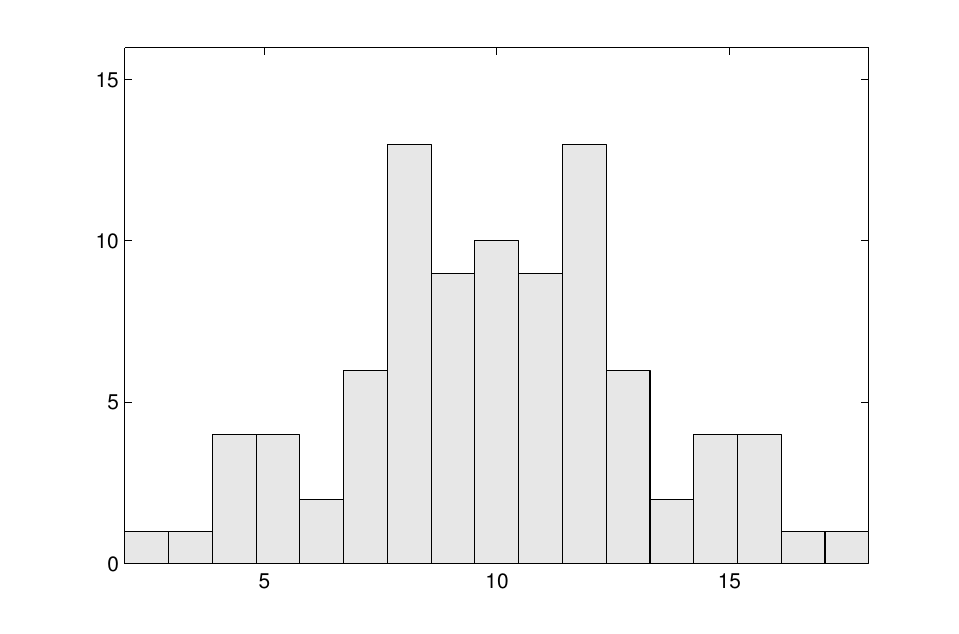}
\\
\includegraphics[scale=.58]{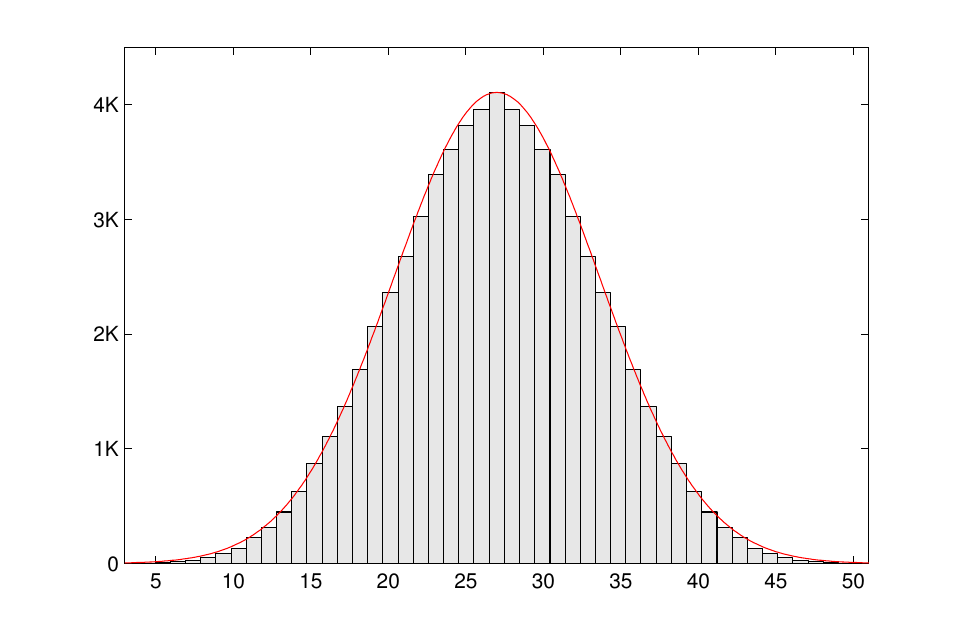}
\end{tabular}
\end{center}
\vspace{-6mm}
\caption{Histograms of the number of inversions in matrices in $\A{4,2}$ (top) and $\A{6,2}$ (bottom).\label{fig:distribution}}
\end{figure}

Figure~\ref{fig:distribution} presents two histograms of values of $\nu(A)$
where $A\in\A{n,2}$, for $n=4$ and $n=6$.
The plot for $n=6$ shows a bell shape (a bell curve is superimposed)
and a similar pattern is observed for $5\le n\le9$).
Therefore, the maximum size of the antichains $\nu^{-1}(t)$ of
Lemma~\ref{lem:inv_antichain} is much larger than the average size
given in Corollary~\ref{cor:average}.
For $n=6$, this maximum value is $|\nu^{-1}(27)|=4108$,
while the average value is approximately $1387$.
Table~\ref{tbl:nu} presents similar computational results for
more values of~$n$.

\begin{table}
\begin{center}
\begin{tabular}{c@{\hspace{3em}}c}
\hline
$n$ & Lower bound for $\w(n,2)$
\\
\hline
3 & 2 \\
4 & 13 \\
5 & 161 \\
6 & 4086 \\
7 & 142468 \\
8 & 7033816 \\
9 & 450066504 \\
\hline
\end{tabular}
\end{center}
\caption{Computational lower bounds for the width $\w(n,2)$ of the Bruhat order of $\A{n,2}$,
obtained by taking the largest antichain of all matrices with equal number of
inversions.\label{tbl:nu}}
\end{table}

It is easy to see that $\w(3,2)=|\nu^{-1}(3)|=|\nu^{-1}(6)|=2$,
and an exhaustive computer search verifies that
$\w(4,2)=|\nu^{-1}(8)|=|\nu^{-1}(12)|=13$.
In the remainder of this section we present constructions of
antichains in the Bruhat order of $\A{n,2}$,
consisting of matrices with equal number of inversions.
Given a matrix $A$, we refer to the matrix obtained by
flipping $A$ in the right/left direction as the {\em conjugate} of $A$
and we denote it by~$\cnj{A}$. Namely, if $A=[a_{ij}]$ is $m$ by $n$, then
$\cnj{A}=[b_{ij}]$ is $m$ by $n$ with $b_{ij}=a_{i,n-j+1}$
for all $1\le i\le m$ and $1\le j\le n$.

The following lemma is an elementary application of  the inclusion--exclusion principle.

\begin{lem}
\label{lem:nuconj}
Let $R=(r_1,\ldots,r_m)$ and $S=(s_1,\ldots,s_n)$ be sequences of
nonnegative integers.
Then for all $A\in\A{R,S}$,\vspace{-1ex}
\[\nu(A)+\nu(\cnj{A})={r_1+\cdots+r_m\choose2}
-\sum_{i=1}^m{r_i\choose2}-\sum_{j=1}^n{s_j\choose2}.\]
\end{lem}


If $A\in\A{n,2}$ such that $A=\cnj{A}$, then by Lemma~\ref{lem:nuconj},
$2\nu(A)={2n\choose2}-2n$, which gives $\nu(A)=n^2-3n/2$.
While this proves that such {\em self-conjugate} matrices do not exist
in $\A{n,2}$ when $n$ is odd,
for even $n$ they give way to constructions of antichains.
Let $n\ge2$ be even and let $A\in\A{n,2}$ with $A=\cnj{A}$.
Let $C$ be the submatrix of $A$ consisting of its first $n/2$ columns.
Then the last $n/2$ columns of $A$ constitute a submatrix equal to $\cnj{C}$. Hence each row of $C$ sums to~$1$ while every column of $C$ sums to~$2$.
We conclude that, the antichain of self-conjugate matrices in $\A{n,2}$
may be identified with the class $\A{R,S}$ where $R$ is the
$n$--vector of all $1$s and $S$ is the $n/2$--vector of all $2$s.
It is known~\cite{GaoMatheis} that the above class $\A{R,S}$ has
size $n!/2^{n/2}$
(see the sequence A000680~\cite{oeis:A000680} in the On-Line Encyclopedia of Integer Sequences).

\begin{thm}
\label{thm:AcEven}
If $n\ge2$ is an even integer, then there is an antichain of size
$\displaystyle\frac{n!}{2^{n/2}}$ in the Bruhat order of $\A{n,2}$.
\end{thm}

\begin{proof}
Let $n=2k$ and let $C$ be a $2k$ by $k$ binary matrix whose
row sums are all $1$, and whose column sums are all $2$.
Let $A_C=\left[C\ \cnj{C}\right]$. Then $A\in\A{n,2}$ and
$\cnj{A}_C=A_C$, thus $\nu(A_C)=n^2-3n/2$.
By Lemma~\ref{lem:inv_antichain}, the set of all these matrices
$A_C$ is an antichain in the Bruhat order of $\A{n,2}$.
This antichain has the desired size since there are precisely
$(2k)!/2^k$ such matrices~$C$.
\end{proof}

A similar construction to that of Theorem~\ref{thm:AcEven}
may be used to obtain an antichain in the Bruhat order of
$\A{n,2}$ when $n$ is odd.

\begin{thm}
If $n\ge3$ is an odd integer, then there is an antichain of size
\label{thm:AcOdd}
$\displaystyle\frac{(n-1)!}{2^{(n-3)/2}}$ in the Bruhat order of $\A{n,2}$.
\end{thm}

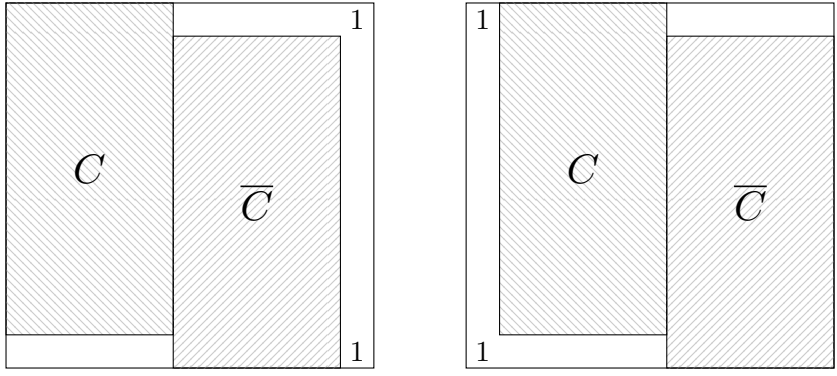
\begin{figure}
\begin{center}
\begin{tabular}{c@{\hspace{12mm}}c}
\begin{tikzpicture}[scale=.44]
\draw (0,0) rectangle (11,11);
\draw[pattern=north west lines,pattern color=gray!40] (0,1) rectangle (5,11);
\draw[pattern=north east lines,pattern color=gray!40] (5,0) rectangle (10,10);
\node[scale=1] at (10.5,10.5) {$1$};
\node[scale=1] at (10.5,.5) {$1$};
\node[scale=1.42] at (2.5,6) {$C$};
\node[scale=1.42] at (7.5,5) {$\cnj{C}$};
\end{tikzpicture}
&
\begin{tikzpicture}[scale=.44]
\draw (0,0) rectangle (11,11);
\draw[pattern=north west lines,pattern color=gray!40] (1,1) rectangle (6,11);
\draw[pattern=north east lines,pattern color=gray!40] (6,0) rectangle (11,10);
\node[scale=1] at (0.5,10.5) {$1$};
\node[scale=1] at (0.5,.5) {$1$};
\node[scale=1.42] at (3.5,6) {$C$};
\node[scale=1.42] at (8.5,5) {$\cnj{C}$};
\end{tikzpicture}
\end{tabular}
\end{center}
\caption{The two matrices $A_C$ and $A_C'$ constructed from $C$
in the proof of Theorem~\ref{thm:AcOdd}.\label{fig:Ac}}
\end{figure}

\begin{proof}
Let $n=2k+1$ and let $C$ be a $2k$ by $k$ binary matrix whose
row sums are all $1$, and whose column sums are all $2$.
Let $A_C=[a_{ij}]$ be the $n$ by $n$ matrix where
\[\begin{aligned}
&A_C[\{1,\ldots,2k\},\{1,\ldots,k\}]=C,\\
&A_C[\{2,\ldots,2k+1\},\{k+1,\ldots,2k\}]=\cnj{C},\\
&a_{1n}=a_{nn}=1,\\
\end{aligned}
\]
and where every entry not included above equals~$0$
(see Figure~\ref{fig:Ac}).
Then $A_C\in\A{n,2}$ and it has 
$\nu(C)+\nu(\cnj{C})={2k\choose2}-k$ inversions within $C$ and $\cnj{C}$,
as well as $2k-1\choose2$ inversions involving one entry from $C$
and one from $\cnj{C}$,
and $4k-1$ inversions involving the $(1,n)$--entry.
Thus we obtain
\[\nu(A_C)={2k\choose2}-k+{2k-1\choose2}+4k-1=k(4k-1).\]
The matrix $A_C'$ obtained from $A_C$ by moving its last column to the
left-most position, has the same number of inversions (counted in a similar
fashion). Since the matrices $A_C$ and $A_C'$ all have the same number of inversions, they
form an antichain in the Bruhat order of $\A{n,2}$.
There are $(2k)!/2^{k}$ binary matrices $C$
of size $2k$ by $k$
with row sums $1$ and column sums $2$, and each
such matrix $C$ contributes two
matrices $A_C$ and $A_C'$ to this antichain.
Therefore, the constructed antichain has size $(2k)!/2^{k-1}$.
\end{proof}

\section{A product construction}

In this section we give a construction of antichains in the Bruhat order
of $\A{R,S}$ that are products of known antichains.
We apply this construction to obtain antichains in the Bruhat
order of $\A{2k,k}$ from the antichains of Section~\ref{sec:An2}.
In the following theorem, $U\otimes V$ denotes the Kronecker product
of the vectors $U$ and $V$. 
For a vector $U=\langle u_1,u_2,\ldots,u_n\rangle$ and a constant
$t\in\R$, the translation vector $\langle t+u_1,t+u_2,\ldots,t+u_n\rangle$
is denoted by $t+U$.

\begin{thm}
For positive integers $a,b,m,n$, let
$R_1=\langle r_1,r_2,\ldots,r_a\rangle$,
$R_2=\langle r'_1,r'_2,\ldots,r'_m\rangle$,
$R_3=\langle r''_1,r''_2,\ldots,r''_m\rangle$,
$S_1=\langle s_1,s_2,\ldots,s_b\rangle$,
$S_2=\langle s'_1,s'_2,\ldots,s'_n\rangle$ and
$S_3=\langle s''_1,s''_2,\ldots,s''_n\rangle$
be vectors of nonnegative integers.
Let $u=r_1+r_2+\cdots+r_a$,
$u'=r'_1+r'_2+\cdots+r'_m$, and $u''=r''_1+r''_2+\cdots+r''_m$,
and suppose that $u'\not=u''$.
If $\mathcal{D}_1$, $\mathcal{D}_2$ and $\mathcal{D}_3$ are antichains
in the Bruhat order of the classes $\A{R_1,S_1}$, $\A{R_2,S_2}$ and
$\A{R_3,S_3}$ respectively,
then there is an antichain of size 
$|\mathcal{D}_1||\mathcal{D}_2|^u|\mathcal{D}_3|^{ab-u}$
in the Bruhat order of the class $\A{R,S}$ where
$R=R_1\otimes R_2+(b-R_1)\otimes R_3$
and $S=S_1\otimes S_2+(a-S_1)\otimes S_3$.
\label{thm:prod}
\end{thm}

\begin{proof}
Let $\mathcal{D}$ be the class of all $am$ by $bn$ binary matrices
$X$ constructed as follows.
For every $P=[p_{ij}]\in\mathcal{D}_1$, we define $X$ as a block matrix
\[
X=\left[\begin{array}{cccc}
X_{11}&X_{12}&\cdots&X_{1b}\\
X_{21}&X_{22}&\cdots&X_{2b}\\
\vdots&\vdots&\ddots&\vdots\\
X_{a1}&X_{a2}&\cdots&X_{ab}\\	
\end{array}\right],
\]
where $X_{ij}\in\mathcal{D}_2$ if $p_{ij}=1$, and
$X_{ij}\in\mathcal{D}_3$ if $p_{ij}=0$.
Clearly,
$|\mathcal{D}|=|\mathcal{D}_1||\mathcal{D}_2|^u|\mathcal{D}_3|^{ab-u}$.
On the other hand, in the $a$ rows of $X$ corresponding to the
$i$--th row of $P$, there are $r_i$ blocks from $\mathcal{D}_2$
and $b-r_i$ blocks from $\mathcal{D}_3$, giving
$r_ir'_k+(b-r_i)r''_k$ ones in the $k$--th of these rows.
Similarly, there are $s_js'_\ell+(a-s_j)s''_\ell$ ones in
the $(j-1)b+\ell$--th column of $X$.
Therefore, $\mathcal{D}\subset\A{R,S}$.
In the rest of the proof we show that every two matrices in $\mathcal{D}$
are incomparable in the Bruhat order of $\A{R,S}$.
For $X,Y\in\mathcal{D}$, let $\Psi_{ij}$ (resp. $\Phi_{ij}$)
be the submatrix of $\Sigma_X$ (resp. $\Sigma_Y$) induced by
rows $(i-1)a+1,(i-1)a+2,\ldots,(i-1)a+a$, and
columns $(j-1)b+1,(j-1)b+2,\ldots,(j-1)b+b$.
We consider two cases.

\emph{Case 1}: The matrices $X,Y\in\mathcal{D}$ correspond to different underlying
matrices $P$ and $Q$ respectively.
For every $1\le k\le a$ and $1\le \ell\le b$, let
$\psi_{k\ell}$ be the bottom-right entry in $\Psi_{k\ell}$.
Note that every block $X_{ij}$ with $i\le k$ and $j\le\ell$
contributes $u'$ to $\psi_{k\ell}$ if $X_{ij}\in\mathcal{D}_2$, 
and it contributes $u''$ to $\psi_{k\ell}$ if $X_{ij}\in\mathcal{D}_3$.
Therefore,
\[\psi_{k\ell}=u'\,\sigma_{k\ell}(P)+u''\big(k\ell-\sigma_{k\ell}(P)\big)
=u''k\ell+(u'-u'')\,\sigma_{k\ell}(P).\]
Similarly, if $\varphi_{k\ell}$ is the bottom-right entry in $\Phi_{k\ell}$,
then
\[\varphi_{k\ell}=u''k\ell+(u'-u'')\,\sigma_{k\ell}(Q).\]
Since $u'-u''\not=0$ and $P$ and $Q$ are incomparable in the Bruhat
order of $\A{R_1,S_1}$, we conclude that $X$ and $Y$ are incomparable
in the Bruhat order of $\A{R,S}$.

\emph{Case 2}: The matrices $X,Y\in\mathcal{D}$ correspond to the same
underlying matrix $P$. Let $1\le i\le a$, $1\le j\le b$,
$1\le k\le m$, $1\le\ell\le n$, $\alpha=a(i-1)+k$, and
$\beta=b(j-1)+\ell$. From the block structure of $X$
we have
\[\begin{aligned}
\sigma_{\alpha\beta}(X) = \sigma_{k\ell}(X_{ij})
   & +\sigma_{i-1,j-1}(P)u'+\big[(i-1)(j-1)-\sigma_{i-1,j-1}(P)\big]u''\\
   & +\big[\sigma_{i,j-1}(P)-\sigma_{i-1,j-1}(P)\big]u'
     +\big[j-1-\sigma_{i,j-1}(P)+\sigma_{i-1,j-1}(P)\big]u''\\
   & +\big[\sigma_{i-1,j}(P)-\sigma_{i-1,j-1}(P)\big]u'
     +\big[i-1-\sigma_{i-1,j}(P)+\sigma_{i-1,j-1}(P)\big]u''.
\end{aligned}\]
All terms on the right hand side of the above equation, except for the
first term, are independent of the individual blocks of~$X$.
A similar formula can be obtained for the matrix $Y$ and putting these
two formulae together, we conclude that
\begin{equation}
\sigma_{\alpha\beta}(X) - \sigma_{\alpha\beta}(Y)
=\sigma_{k\ell}(X_{ij}) - \sigma_{k\ell}(Y_{ij}).
\label{eq:sigma}
\end{equation}
Since $X$ and $Y$ correspond to the same matrix $P\in\mathcal{D}_1$,
for each $1\le i\le a$ and $1\le j\le b$, the blocks $X_{ij}$
and $Y_{ij}$ belong to the same antichain $\mathcal{D}_2$ or
$\mathcal{D}_3$. If $X\not=Y$, then for some $1\le i\le a$ and
$1\le j\le b$, $X_{ij}$ and $Y_{ij}$ are inequal, thus incomparable.
This implies incomparability of $X$ and $Y$ by equation~(\ref{eq:sigma}).
\end{proof}

\begin{cor}
Let $a$ and $m$ be positive integers where $a$ is even and let $n=am$.
Let $\mathcal{D}_1$ be an antichain in the Bruhat order of $\A{a,a/2}$
and $\mathcal{D}_2$ be an antichain in the Bruhat order of $\A{m,2}$.
Then there is an antichain of size $|\mathcal{D}_1||\mathcal{D}_2|^{a^2}$
in the Bruhat order of $\A{n,n/2}$.
\label{cor:2kkproduct}
\end{cor}

\begin{proof}
Let 
$\mathcal{D}_3=\big\{J_m-X\,\big|\ X\in\mathcal{D}_2\big\}$
be the class of matrix complements of members of $\mathcal{D}_2$.
It is easy to see that matrix complement reverses the Bruhat order,
thus it preserves incomparability.
Therefore, $\mathcal{D}_3$ is an antichain in the Bruhat order of
$\A{n,n-2}$. The desired antichain is obtained by applying
Theorem~\ref{thm:prod} to the antichains
$\mathcal{D}_1$, $\mathcal{D}_2$ and $\mathcal{D}_3$.
\end{proof}

Conflitti\etal~\cite{ConflittiDaFonsecaMamede2013}
give constructions of antichains of size $(\floor{k/2}^4+1)^2$
in the Bruhat order of $\A{2k,k}$. The following improves this result.

\begin{cor}
\label{cor:2kk}
Let $k$ be a positive integer and let $\w(2k,k)$
denote the width of the Bruhat order of $\A{2k,k}$. Then
\[\w(2k,k)\ge g(k)=\begin{cases}
\displaystyle\frac{(k!)^4}{{4^k}} & \text{\ if }k\text{ is even,}\\[1em]
\displaystyle\frac{\big[(k-1)!\big]^4}{{4^{k-3}}} & \text{\ if }k\text{ is odd.}\\
\end{cases}\]
\end{cor}

\begin{proof}
An antichain of size $g(k)$ in the Bruhat order of $\A{2k,k}$
is obtained by applying Corollary~\ref{cor:2kkproduct} to
the antichains $\mathcal{D}_1=\{I_2\}$ in the Bruhat order of
$\A{2,1}$, and  $\mathcal{D}_2$ the Bruhat order of $\A{k,2}$
given in Theorem~\ref{thm:AcEven} or~\ref{thm:AcOdd},
depending on the parity of~$k$.
\end{proof}

\begin{rem}
The lower bound of Corollary~\ref{cor:2kk} may be improved by
using other antichains $\mathcal{D}_1$. For example,
if $k$ is a multiple of $4$, we may take $\mathcal{D}_1$ to be
an antichain of size $13$ in the Bruhat order of $\A{4,2}$
and $\mathcal{D}_2$ the antichain of Theorem~\ref{thm:AcEven}
in the Bruhat order of $\A{k/2,2}$ to obtain an antichain of
size $13[(k/2)!]^{16}/16^k$ in the Bruhat order of $\A{2k,k}$.
\label{rem:improve}
\end{rem}

\section{Concluding remarks}

We give constructions of antichains in the Bruhat order of the
classes $\A{n,2}$ and $\A{2k,k}$, where $n$ and $k$ are positive
integers. These antichains give asymptotically exponential lower
bounds for the width of the Bruhat order of these classes.
Remark~\ref{rem:improve} and the discussion following
Corollary~\ref{cor:average} indicate that these bounds can be
further improved.
The main tool in our constructions of antichains in the Bruhat
order of $\A{n,2}$ is the fact that if two matrices in this class
have the same number of inversions, then they are incomparable in the
Bruhat order. An affirmative answer to the following problem
would provide the same tool for more classes $\A{R,S}$.
This problem is a weakened restatement of Question~4
in~\cite{Mamad:An2Chain}.

\begin{pr}
If $A,C\in\A{R,S}$, does $\nu(A)=\nu(C)$ imply that $A$ and
$C$ are incomparable in the Bruhat order of $\A{R,S}$?
\end{pr}


\def\soft#1{\leavevmode\setbox0=\hbox{h}\dimen7=\ht0\advance \dimen7
  by-1ex\relax\if t#1\relax\rlap{\raise.6\dimen7
  \hbox{\kern.3ex\char'47}}#1\relax\else\if T#1\relax
  \rlap{\raise.5\dimen7\hbox{\kern1.3ex\char'47}}#1\relax \else\if
  d#1\relax\rlap{\raise.5\dimen7\hbox{\kern.9ex \char'47}}#1\relax\else\if
  D#1\relax\rlap{\raise.5\dimen7 \hbox{\kern1.4ex\char'47}}#1\relax\else\if
  l#1\relax \rlap{\raise.5\dimen7\hbox{\kern.4ex\char'47}}#1\relax \else\if
  L#1\relax\rlap{\raise.5\dimen7\hbox{\kern.7ex
  \char'47}}#1\relax\else\message{accent \string\soft \space #1 not
  defined!}#1\relax\fi\fi\fi\fi\fi\fi}

\end{document}